\documentclass[11pt]{amsart}
\usepackage{amssymb, amsmath, amsthm, amscd}
\newtheorem{proposition}{Proposition}[section]

\newtheorem{theorem}[proposition]{Theorem}

\theoremstyle{definition}

\newtheorem{remark}[proposition]{Remark}

\newcommand{\thlabel}[1]{\label{th:#1}}
\newcommand{\thref}[1]{Theorem~\ref{th:#1}}
\newcommand{\selabel}[1]{\label{se:#1}}

\newcommand{\relabel}[1]{\label{re:#1}}

\newcommand{\eqlabel}[1]{\label{eq:#1}}
\newcommand{\equref}[1]{(\ref{eq:#1})}

\def\ot{\otimes}

\newcommand{\Cc}{\mathcal{C}}

\def\*C{{}^*\hspace*{-1pt}{\Cc}}
\def\text#1{{\rm {\rm #1}}}

\input xy
\xyoption {all} \CompileMatrices

\begin{document}

\title[]{Categorical constructions for Hopf algebras}
\dedicatory{Dedicated to the memory of Professor Liliana Pavel}

\author{A.L. Agore}\thanks{The author acknowledges partial support
from CNCSIS grant 24/28.09.07 of PN II "Groups, quantum groups,
corings and representation theory". }
\address{Department of Mathematics, Academy of Economic Studies, Piata
Romana 6, RO-010374 Bucharest 1,  Romania}
\email{ana.agore@fmi.unibuc.ro}



\begin{abstract}
We prove that both, the embedding of the category of Hopf algebras
into that of bialgebras and the forgetful functor from the
category of Hopf algebras to the category of algebras, have right
adjoints; in other words: every bialgebra has a Hopf coreflection
and on every algebra there exists a cofree Hopf algebra. In this
way we give an affirmative answer to a forty years old problem
posed by Sweedler. On the route the coequalizers and the
coproducts in the category of Hopf algebras are explicitly
described.
\end{abstract}

\maketitle

\section*{Introduction}

Hopf algebras appeared naturally in the study of Lie groups
cohomology. The survey paper \cite{AS} covers the beginnings of
Hopf algebras and the roles played by H. Hopf, P. Cartier, A.
Borel, J. Milnor, J. Moore, B. Konstant and M. Sweedler in the
development of this theory. Hopf algebras became a fervid field of
study especially after the appearance of the monograph \cite{Sw}.
In the present paper we bring new contributions to the study of
the category of Hopf algebras.

We turn our attention to the fundamental book of Sweedler: in
\cite[p. 135]{Sw} are stated, without any proofs, the following
problems concerning Hopf algebras: given a coalgebra $C$ there
exists a free Hopf algebra on $C$ (i.e. the forgetful functor from
the category of Hopf algebras to the category of coalgebras has a
left adjoint) and a free commutative Hopf algebra on $C$. The
problem has turned out to be quite difficult: several years passed
until Takeuchi, in \cite[Sections \S1 and \S11]{T}, answered
affirmatively to both statements. His proof relies on an ingenious
and laborious construction. Moreover, in \cite[p. 135]{Sw}
Sweedler also states, again without any proofs, the dual of the
above problem: given an algebra $A$ there exists a cofree Hopf
algebra on $A$ (that is the forgetful functor from the category of
Hopf algebras to the category of algebras has a right adjoint) and
a cofree cocommutative Hopf algebra on $A$. Concerning this
problem, recently H.-E. Porst \cite[Corollary 4.1.4]{HP} proved
that the existence of a cofree Hopf algebra on every algebra
implies the existence of a cofree cocommutative Hopf algebra on
every algebra. In the present paper we prove Sweedler's statement
concerning the existence of a cofree Hopf algebra on every
algebra.

The paper is structured as follows. In Section $1$ we introduce
the notations and recall, without proofs, some well known results
pertaining to category theory that will be intensively used
throughout the paper. In Section $2$ we give an explicit
description of coequalizers and coproducts in the category
$k$-HopfAlg of Hopf algebras. Using the aforementioned
constructions we prove, using the Special Adjoint Functor Theorem,
Sweedler's statement concerning the existence of a cofree Hopf
algebra on every algebra (\thref{S}).

\section{Preliminaries}\selabel{1}

Throughout this paper, $k$ will be a field. Unless specified
otherwise, all vector spaces, algebras, coalgebras, bialgebras,
tensor products and homomorphisms are over $k$. Our notation for
the standard categories is as follows: ${}_k{\mathcal {M}}$
($k$-vector spaces), $k$-Alg (associative unital $k$-algebras),
$k$-BiAlg (bialgebras over $k$), $k$-HopfAlg (Hopf algebras over
$k$). We refer to \cite{Sw} for further details concerning Hopf
algebras.

We use the standard notations for opposite and coopposite
structures: $A^{op}$ denotes the opposite of the algebra $A$ and
$C^{cop}$ stands for the coopposite of the coalgebra $C$.

Let us recall briefly some well known results from category
theory, refering the reader to \cite{ML} for more details. A
category $\mathcal{C}$ is called \textit{(co)complete} if all
diagrams in $\mathcal{C}$ have (co)limits in $\mathcal{C}$. A
category $\mathcal{C}$ is \textit{(co)complete} if and only if
$\mathcal{C}$ has (co)equalizers of all pairs of arrows and all
(co)products \cite[Theorem 6.10]{par}. The category $\mathcal{C}$
is called \textit{locally small} (or \textit{well-powered}) if the
subobjects of each $C \in \mathcal{C}$ can be indexed by a set.
Dually, the category $\mathcal{C}$ is \textit{colocally small} (or
\textit{co-well-powered}) if its dual is locally small. At some
point we will also use, in passing, the notion of locally
presentable category. More details regarding this type of
categories can be found in \cite{adamek}. All categories
considered above are locally presentable. Thus, they are
cocomplete by the definition of locally presentable categories and
complete by \cite[Remark 1.56]{adamek}. For a more detailed
discussion concerning the completeness and cocompleteness of the
above categories we refer the reader to \cite{HP}, \cite{HP1} and
\cite{HP2}. A subcategory $\mathcal{D}$ of $\mathcal{C}$ is called
\textit{(co)reflective} in $\mathcal{C}$ when the inclusion
functor $U: \mathcal{D} \rightarrow \mathcal{C}$ has a (right)left
adjoint.

The following categorical result play a key role in showing that
the category of Hopf algebras is a coreflective subcategory of the
category of bialgebras:

\begin{theorem}(The Special Adjoint Functor Theorem)\thlabel{02}
If $\mathcal{C}$ is a complete and locally small category with a
cogenerator, then a functor $G: \mathcal{C} \rightarrow
\mathcal{D}$ has a left adjoint if and only if it is limit
preserving. Dually, if $\mathcal{C}$ is a cocomplete and colocally
small category with a generator, then a functor $G: \mathcal{C}
\rightarrow \mathcal{D}$ has a right adjoint if and only if it is
colimit preserving.
\end{theorem}

\section{Cofree Hopf algebras generated by algebras}\selabel{2}

Recall that the forgetful functor from the category of groups to
the category of monoids has a left adjoint, the so-called
anvelopant group of a monoid, and a right adjoint, which assigns
to each monoid the group of its invertible elements. Therefore, if
we think of Hopf algebras as a natural generalization of groups,
we may expect the same behavior in the case of the embedding
functor $F$: $k$-HopfAlg $\rightarrow$ $k$-BiAlg from the category
of Hopf algebras to the category of bialgebras. It is well known
that the above embedding functor $F$ has a left adjoint
\cite[Theorem 2.6.3]{P}. We shall prove in this section that $F$
has also a right adjoint. This result together with that fact that
there exist a cofree bialgebra on every algebra led us to the
conclusion that the forgetful functor $F$: $k$-HopfAlg
$\rightarrow$ $k$-Alg has a right adjoint, as stated in \cite{Sw}.
In order to prove our main result we need some preparations. We
start by recalling, for a further use, the constructions of
coequalizers and coproducts in the category $k$-BiAlg of
bialgebras.

Let $(A,m_{A},\eta_{A},\Delta_{A},\varepsilon_{A})$,
$(B,m_{B},\eta_{B},\Delta_{B},\varepsilon_{B})$ be two bialgebras
and $f$, $g: B \rightarrow A$ be two bialgebra maps. Consider $I$
the two-sided ideal generated by $\{f(b)-g(b) \, | \, b \in B \,
\}$. By a simple computation it can be seen than $I$ is also a
coideal. Then $(A/I, \pi)$ is the coequalizer of the morphisms
$(f,g)$ in $k$-BiAlg, where $\pi:A \rightarrow A/I$ is the
canonical projection. Indeed for all bialgebras $H$ and all
bialgebra morphisms $h: A \rightarrow H$ such that $h \circ f = h
\circ g$ we obtain $I \subseteq kerh$, hence there exists an
unique bialgebra map $h' : A/I \rightarrow H$ such that $h' \circ
\pi = h$.

\begin{remark}\relabel{1}
Note that if $A$, $B$ are two Hopf algebras and $f, g: B
\rightarrow A$ are Hopf algebra maps then the ideal $I$ defined
above is actually  a Hopf ideal and $(A/I, \pi)$ is the
coequalizer of the morphisms $(f,g)$ in $k$-HopfAlg, where $\pi:A
\rightarrow A/I$ is the canonical projection.
\end{remark}

Next, we recall from \cite{P} the construction of the coproduct in
the category $k$-BiAlg of bialgebras. Let $(A_{l})_{l \in I}$ be a
family of algebras, $\bigl(\bigoplus_{l \in I}A_{l}, (j_{l})_{l
\in I}\bigl)$ be the coproduct in ${}_k{\mathcal {M}}$ and $i:
\bigoplus_{l \in I}A_{l} \rightarrow T\bigl(\bigoplus_{l \in
I}A_{l}\bigl)$ be the canonical inclusion where $T
\bigl(\bigoplus_{l \in I}A_{l}\bigl)$ is the tensor algebra of the
vector space $\bigoplus_{l \in I}A_{l}$. Then $\bigl(\coprod_{l
\in I} A_{l} := T\bigl(\bigoplus_{l \in I}A_{l}\bigl)/L ,
(q_{l})_{l \in I}\bigl)$ is the coproduct of the above family in
$k$-Alg, where $L$ is the two sided ideal in $T\bigl(\bigoplus_{l
\in I}A_{l}\bigl)$ generated by the set $J:= \{i\circ
j_{l}(x_{l}y_{l})-i\bigl(j_{l}(x_{l})\bigl)i\bigl(j_{l}(y_{l})\bigl)
,$ $1_{T\bigl(\bigoplus A_{l}\bigl)} - i\circ j_{l}(1_{A_{l}}) |
x_{l}, y_{l} \in A_{l}, l \in I \}$, $\nu : T\bigl(\bigoplus_{l
\in I}A_{l}\bigl) \rightarrow T\bigl(\bigoplus_{l \in
I}A_{l}\bigl)/L$ denotes the canonical projection and $q_{l} = \nu
\circ i \circ j_{l}$ for all $l \in I$. Furthermore, $\coprod_{l
\in I} A_{l}$ is actually a bialgebra provided that $(A_{l})_{l
\in I}$ is a family of bialgebras. The comultiplication and the
counit are given by the unique algebra maps such that the
following diagrams commute:
\begin{equation}\eqlabel{B}
\xymatrix {& A_{l} \ar[r]^{q_{l}} \ar[dr]_{(q_{l} \otimes q_{l})\circ \Delta_{l}} & {\coprod_{l \in I} A_{l}} \ar[d]^{\Delta}\\
& {} & {\coprod_{l \in I} A_{l} \otimes \coprod_{l \in I} A_{}}}
\xymatrix {& A_{l} \ar[r]^{q_{l}} \ar[dr]_{\varepsilon_{l}} & {\coprod_{l \in I} A_{l}} \ar[d]^{\varepsilon}\\
& {} & k }
\end{equation}
and $\coprod_{l \in I} A_{l}$ is the coproduct in $k$-BiAlg of the
above family of bialgebras (\cite[Corollary 2.6.2]{P}).

Now let $\bigl(H_{l}, m_{l}, \eta_{l}, \Delta_{l},
\varepsilon_{l}, S_{l}\bigl)_{l \in I}$ be a family of Hopf
algebras. Consider $\bigl((H:=\coprod_{l \in I} H_{l},$ $m, \eta,
\Delta, \varepsilon), (q_{l})_{l \in I}\bigl)$ the coproduct of
the above family in the category $k$-BiAlg of bialgebras.
\newline The universal property of the coproduct yields
an unique bialgebra map $S:H \rightarrow H^{opcop}$ such that the
following diagram commutes for all $l \in I$:
\begin{equation}\eqlabel{antipod}
\xymatrix {& H_{l} \ar[r]^{q_{l}} \ar[d]_{S_{l}} & {H} \ar[d]^{S}\\
& {{H_{l}}^{opcop}}\ar[r]^{q_{l}} & {H^{opcop}}}
\end{equation}

With the notations above we have the following result which
provides a completely description of the coproducts in the
category $k$-HopfAlg of Hopf algebras:

\begin{theorem}\thlabel{4}
Let $\bigl(H_{l}, m_{l}, \eta_{l}, \Delta_{l}, \varepsilon_{l},
S_{l}\bigl)_{l \in I}$ be a family of Hopf algebras. The Hopf
algebra $\bigl(H:=\coprod_{l \in I} H_{l},$ $m, \eta, \Delta,
\varepsilon, S\bigl)$ together with structure maps $(q_{l})_{l \in
I}$ is the coproduct in the category $k$-HopfAlg of the family
$\bigl(H_{l}, m_{l}, \eta_{l}, \Delta_{l}, \varepsilon_{l},
S_{l}\bigl)_{l \in I}$ of Hopf algebras. In particular, the
category $k$-HopfAlg is cocomplete.
\end{theorem}

\begin{proof}
We will first prove that $S$ is an antipode for the bialgebra $H$,
i.e.
\begin{equation}\eqlabel{S}
m \circ \bigl(Id \ot S\bigl) \circ \Delta = m \circ \bigl(S \ot
Id\bigl) \circ \Delta = \eta \circ \varepsilon
\end{equation}
Since $S:H \rightarrow H^{opcop}$ defined in \equref{antipod} is a
bialgebra map we only need to prove that \equref{S} holds only on
the generators of $H$ as an algebra. Indeed, let $h, k$ be
generators in $H$ for which \equref{S} holds. We obtain :
$$(hk)_{(1)}S((hk)_{(2)}) = h_{(1)}k_{(1)}S(k_{(2)})S(h_{(2)}) =
\varepsilon(k) h_{(1)}S(h_{(2)}) =
\varepsilon(h)\varepsilon(k)1_{H} = \varepsilon(hk)1_{H}$$ It
follows from here that \equref{S} also holds for $kh$ and thus it
holds for all elements
in $H$.\\
Now having in mind that $H := T\bigl(\bigoplus_{l \in
I}H_{l}\bigl)/L$ we only need to prove \equref{S} for the elements
$\widehat{x} \in H$ with $x \in \bigoplus_{l \in I}H_{l}$, whereas
the tensor algebra $T\bigl(\bigoplus_{l \in I}H_{l}\bigl)$ is the
free algebra on $\bigoplus_{l \in I}H_{l}$. Moreover, since
$\bigoplus_{l \in I}H_{l} = \{x \in \prod_{l \in I}H_{l} ~|~
{\rm~supp~}(x) < \infty\}$ it is enough to show that \equref{S}
holds for all $x_{l} \in H_{l}$, $l \in I$. We then have:
\begin{eqnarray*}
m \circ \bigl(Id \ot S\bigl) \circ \Delta (\widehat{x_{l}}) &{=}&
m \circ \bigl(Id \ot S\bigl) \circ \Delta \circ q_{l}(x_{l})\\
&\stackrel{\equref{B}} {=}& m \circ \bigl(Id \ot S\bigl) \circ
(q_{l} \otimes q_{l})\circ \Delta_{l}(x_{l})\\
&{=}& m \circ \bigl(q_{l} \ot (S \circ q_{l})\bigl)\circ
\Delta_{l}(x_{l})\\
&\stackrel{\equref{antipod}} {=}& m \circ \bigl(q_{l} \ot (q_{l}
\circ
S_{l})\bigl)\circ \Delta_{l}(x_{l})\\
&{=}&m \circ (q_{l} \ot q_{l}) \circ (Id \circ
S_{l})\circ \Delta_{l}(x_{l})\\
&\stackrel{q_{l} - {\rm~algebra~map~}}{=}&q_{l} \circ m_{l} \circ
(Id \circ
S_{l})\circ \Delta_{l}(x_{l})\\
&{=}&q_{l} \circ \eta_{l} \circ \varepsilon_{l}(x_{l})\\
&\stackrel{q_{l} - {\rm~algebra~map~}}{=}&\eta \circ \varepsilon_{l}(x_{l})\\
&\stackrel{q_{l} - {\rm~coalgebra~map~}}{=}&\eta \circ \varepsilon
\circ
q_{l}(x_{l})\\
&{=}&\eta \circ \varepsilon(\widehat{x_{l}})
\end{eqnarray*}
Hence $m \circ \bigl(Id \ot S\bigl) \circ \Delta = \eta \circ
\varepsilon$. In the same way it can be proved that $m \circ
\bigl(S \ot Id\bigl) \circ \Delta = \eta \circ \varepsilon$. Thus
$S$ is an antipode for $H$, as desired.

Now since $k$-HopfAlg is a full subcategory of the category
$k$-BiAlg it follows that $\bigl((H:=\coprod_{l \in I} H_{l},$ $m,
\eta, \Delta, \varepsilon), (q_{l})_{l \in I}\bigl)$ is also the
coproduct of the family $\bigl(H_{l}, m_{l}, \eta_{l}, \Delta_{l},
\varepsilon_{l}, S_{l}\bigl)_{l \in I}$ of Hopf algebras in the
category $k$-HopfAlg.

\end{proof}

We need the following well known result:

\begin{theorem}(\cite[page 134]{Sw})\thlabel{31}
The forgetful functor $F: k-BiAlg \rightarrow k-Alg$ has a right
adjoint, i.e. there exists a cofree bialgebra on every algebra.
\end{theorem}

Our main results now follow:

\begin{theorem}\thlabel{32}
The embedding functor $F: k-HopfAlg \rightarrow k-BiAlg$ has a
right adjoint, i.e. the category of Hopf algebras is a
coreflective subcategory of the category of bialgebras.
\end{theorem}

\begin{proof}
We apply the Special Adjoint Functor Theorem (\thref{02}): since
$k$-HopfAlg is a locally presentable category by \cite[4.3.1 and
4.1.3]{HP} (although epimorphisms of Hopf algebras are not
necessarily surjective maps \cite{CH}), this category in
particular has a generator and is cocomplete and is colocally
small. By \thref{4} the result follows.
\end{proof}

\begin{theorem}\thlabel{S}
The forgetful functor $F: k-HopfAlg \rightarrow k-Alg$ has a right
adjoint, i.e. there exists a cofree Hopf algebra on every algebra.
\end{theorem}

\begin{proof}
It follows from \thref{31} and \thref{32} by composing the right
adjoint functors.
\end{proof}

We proved, using the Special Adjoint Functor Theorem, the
existence of a cofree Hopf algebra on every bialgebra without
indicating explicitly his construction. The following natural
problem arises:

\textbf{Problem:} \textit{Give an explicit construction of the
cofree Hopf algebra on an bialgebra (resp. algebra).}

We expect that the right adjoint of the embedding functor from the
category of Hopf algebras to the category of bialgebras to assign
to every bialgebra $B$ his "biggest" subbialgebra $H$ that has an
antipode.

\section*{Acknowledgements}

The author wishes to thank Professor Gigel Militaru, who suggested
the problem studied here, for his great support and for the useful
comments from which this manuscript has benefitted, as well as the
referee for valuable suggestions and for indicating the papers
\cite{HP1} and \cite{HP2}.

\end{document}